\documentclass[11pt]{article}
\usepackage{amsmath,amssymb} 

\begin{document}

\title{Dehn surgery on knots of wrapping number 2}

\author{Ying-Qing Wu} 
\date{}
\maketitle



\begin{abstract}
  Suppose $K$ is a hyperbolic knot in a solid torus $V$ intersecting a
  meridian disk $D$ twice.  We will show that if $K$ is not the
  Whitehead knot and the frontier of a regular neighborhood of $K \cup
  D$ is incompressible in the knot exterior, then $K$ admits at most
  one exceptional surgery, which must be toroidal.  Embedding $V$ in
  $S^3$ gives infinitely many knots $K_n$ with a slope $r_n$
  corresponding to a slope $r$ of $K$ in $V$.  If $r$ surgery on $K$
  in $V$ is toroidal then either all but at most three $K_n(r_n)$ are
  toroidal, or they are all reducible or small Seifert fibered with
  two common singular fiber indices.  These will be used to classify
  exceptional surgeries on wrapped Montesinos knots in solid torus,
  obtained by connecting the top endpoints of a Montesinos tangle to
  the bottom endpoints by two arcs wrapping around the solid torus.
\end{abstract}

\newcommand{\proof}{\noindent {\bf Proof.} }
\newcommand{\qed}{\quad $\Box$}
\newtheorem{thm}{Theorem}[section]
\newtheorem{prop}[thm]{Proposition} 
\newtheorem{lemma}[thm]{Lemma} 
\newtheorem{sublemma}[thm]{Sublemma} 
\newtheorem{cor}[thm]{Corollary} 
\newtheorem{defn}[thm]{Definition} 
\newtheorem{convention}[thm]{Convention} 
\newtheorem{notation}[thm]{Notation} 
\newtheorem{qtn}[thm]{Question} 
\newtheorem{example}[thm]{Example} 
\newtheorem{remark}[thm]{Remark} 
\newtheorem{conj}[thm]{Conjecture} 
\newtheorem{prob}[thm]{Problem} 
\newtheorem{rem}[thm]{Remark} 

\newcommand{\bdd}{\partial}
\newcommand{\Int}{{\rm Int}}
\newcommand{\wind}{{\rm wind}}
\newcommand{\wrap}{{\rm wrap}}
\renewcommand{\a}{\alpha}
\renewcommand{\b}{\beta}

\input epsf.tex

\section{Introduction} \label{sec:1}

A Dehn surgery on a hyperbolic knot $K$ in a compact 3-manifold is
{\it exceptional\/} if the surgered manifold is non-hyperbolic.  When
the manifold is a solid torus, the surgery is exceptional if and only
if the surgered manifold is either a solid torus, reducible, toroidal,
or a small Seifert fibered manifold whose orbifold is a disk with two
cone points.  Solid torus surgeries have been classified by Berge
\cite{Be} and Gabai \cite{Ga1, Ga2}, and by Scharlemann \cite{Sch}
there is no reducible surgery.  For toroidal surgery, Gordon and
Luecke \cite{GL2} showed that the surgery slope must be either an
integral or a half integral slope.  By \cite{GW1}, this is also true
for small Seifert fibered surgeries.
 
In this paper we study Dehn surgery on hyperbolic knots $K$ in a solid
torus $V$ with wrapping number 2.  The {\em wrapping number\/}
$\wrap(K)$ of a knot $K$ in a solid torus $V$ is defined to be the
minimal geometric intersection number of $K$ with a meridional disk
$D$ of $V$, and the {\em winding number\/} $\wind(K)$ of $K$ is the
algebraic intersection number of $K$ with $D$.  Thus if $K$ is a
knot in a solid torus $V$ with $\wrap(K) = 2$ then $\wind(K) = 0$ or
$2$.  It follows from the results above that there is no reducible or
solid torus surgery on such a hyperbolic knot.  We would like to know
if there are toroidal or small Seifert fibered surgeries on such a
knot.

Exceptional surgery does exist on some knots with wrapping number 2.
If $K$ in $V$ has a spanning surface which is either a punctured torus
or a punctured Klein bottle then surgery on $K$ along the boundary
slope of this surface is toroidal.  A well known example of knots in
solid tori that admit multiple exceptional surgeries is the Whitehead
knot, obtained by deleting an open neighborhood of a component of the
Whitehead link in $S^3$.  It admits a total of 5 exceptional
surgeries, two toroidal and three small Seifert fibered.  

For the case of knots with winding number 2, consider the knot
obtained by putting a Montesinos tangle $T[-1/2, 1/3]$ horizontally in
the solid torus $V$ and then connecting the top endpoints to the
bottom endpoints by two strings running around the solid torus; see
Figure 5.1(b), where $V$ is the complement of the dotted circle.  It
is called a wrapped Montesinos knot and denoted by $K^1(-1/2,\, 1/3)$;
see Section \ref{sec:5} for more details.  We will show that this knot
admits three exceptional surgeries, two toroidal and one small Seifert
fibered.  See Proposition \ref{prop:2.2}.  We suspect that these are
the only examples of knots with wrapping number 2 that admit multiple
exceptional surgeries.

\begin{conj} Suppose $K$ is a hyperbolic knot in a solid torus $V$,
  and $K$ is not the Whitehead knot or the wrapped Montesinos knot
  $K^1(-1/2, 1/3)$.  Then $K$ admits no small Seifert fibered surgery
  and at most one toroidal surgery.
\end{conj}

Let $D$ be a meridian disk intersecting $K$ twice.  Cutting $(V,K)$
along $D$ produces a 2-string tangle $(B, \tau)$.  Let $X$ be the
tangle space $B - \Int N(\tau)$, and let $\bdd_h X$ be the frontier of
$X$ in $V$.  It can be shown that for the knot $K = K^1(1/2,\, 1/3)$
above, this surface $\bdd_h X$ is compressible.  This is a very
special property since most 2-string tangle spaces have incompressible
boundary.  For example, if $\tau$ is a Montesinos tangle of length at
least 2 then $\bdd (B - \Int N(\tau))$ is incompressible unless $\tau$
is equivalent to $T[1/2, p/q]$; see \cite{Wu96}.  The following
theorem shows that the above conjecture is true if $\bdd_h X$ is
incompressible.  Denote by $(V,K,r)$ the manifold obtained by $r$
surgery on a knot $K$ in a 3-manifold $V$.

\bigskip \noindent
{\bf Theorem \ref{thm:3.8}.} 
{\em 
  Suppose $K$ is a hyperbolic knot with $\wrap(K) = 2$ in a solid
  torus $V$, $K$ is not the Whitehead knot, and $\bdd_h X$ is
  incompressible in $X$.  Then $K$ admits at most one exceptional
  surgery $(V,K,r)$, which must be a toroidal surgery and $r$ an
  integral slope.}
\bigskip

Note that the surface $\bdd_h X$ is always incompressible if $K$ is
hyperbolic and $\wind(K) = 0$, hence Conjecture 1.1 is true for 
knots with $\wind(K) = 0$.  Since the Whitehead knot admits 5
exceptional surgeries, it is surprising to see that no other knots
with $\wind(K)=0$ and $\wrap(K)=2$ has more than one exceptional
surgeries.

We now consider knots obtained by embedding $(V,K)$ in the 3-sphere.
Let $\varphi_0$ be a standard embedding, and $\varphi_n$ the
composition of $\varphi_0$ with $n$ right hand full twists of $V$
along a meridian disk.  Denote by $K_n = \varphi_n(K)$ and by $r_n =
\varphi_n(r)$, for a fixed slope $r$ of $K$.  Denote by $K_n(r_n)$ the
surgery on $K_n$ along the slope $r_n$.  Clearly $K_n(r_n)$ is
obtained by Dehn filling $(V,K,r)$ on $\bdd V$, hence if $(V,K,r)$ is
hyperbolic then most $K_n(r_n)$ are hyperbolic.  In general it might
be possible that $(V,K,r)$ is nonhyperbolic while infinitely many
$K_n(r_n)$ are hyperbolic.  However, we will show that this does not
happen when $\wrap(K) = 2$.

\bigskip \noindent
{\bf Theorem \ref{thm:4.3}.} 
{\em  Suppose $\wrap(K) = 2$, and $(V,K,r)$ is non hyperbolic.
  Then $K_n(r_n)$ is nonhyperbolic for all but at most three $n$.
  Moreover, either

  (1) there is an $n_0$ such that $K_n(r_n)$ is toroidal unless $|n -
  n_0| \leq 1$; or

  (2) $K_n(r_n)$ is atoroidal for all $n$, and there exist $q_1, q_2
  \in \Bbb Z$ such each $K_n(r_n)$ is either reducible or has a small
  Seifert fibration with $q_1, q_2$ as the indices of two of its
  singular fibers. }
\bigskip

Thus if $(V,K,r)$ is nonhyperbolic then $K_n(r_n)$ is either toroidal
for all but at most three $n$, or is never toroidal.  This property is
useful in determining whether $(V,K,r)$ is hyperbolic; see the proof
of Theorem 5.5.

Up to homeomorphism there are essentially two ways to make wrapped
Montesinos links from a Montesinos tangle $T[t_1, ..., t_k]$, denoted
by $K^0[t_1, ..., t_k]$ and $K^1[t_1, ..., t_k]$.  See Section
\ref{sec:5} for detailed definitions.  The above theorems will be used
to prove the following classification theorem, which shows that there
is no other exceptional Dehn surgeries on wrapped Montesinos knots in
solid tori besides the well known examples and the ones mentioned
above.  In particular, Conjecture 1.1 is true for these knots.  Here
two pairs $(K,r)$ and $(K',r')$ are equivalent if there is an obvious
homeomorphism of $V$ taking one to the other; see Section 5 for
detailed definitinos.  We may assume that $K \neq K^a[0]$ or
$K^a[1/q]$ as otherwise $K$ is nonhyperbolic.

\bigskip \noindent
{\bf Theorem \ref{thm:5.5}.} 
{\em 
  Suppose $K = K^a(t_1, ..., t_k) \subset V$ is not equivalent to
  $K^a(0)$ or $K^a(1/q)$ for any integer $q$.  Let $(V,K,r)$ be an
  exceptional surgery.  Then $(K, r)$ is equivalent to one of the
  following pairs.  The surgery is small Seifert fibered for $r=1,2,3$
  in (1) and $r=7$ in (4), and toroidal otherwise.

  (1) $K = K^0(2)$ (the Whitehead knot), $r = 0, 1, 2, 3, 4$.

  (2) $K = K^a(n)$, $n>2$, $r=0$ if $a=0$, and $r=2n$ otherwise.

  (3) $K = K^a(1/q_1,\, 1/q_2)$, $|q_i|\geq 2$, and $r$ is the pretzel
  slope.

  (4) $K = K^1(-1/2,\, 1/3)$, $r = 6,7,8$.
}
\bigskip

These results will be used to study Seifert fibered surgery on
Montesinos knots in $S^3$.  We will show that $6+4n$ and $7+4n$
surgeries on hyperbolic $(-2,3,2n+1)$ pretzel knots are Seifert
fibered.  See Corollary \ref{cor:2.3} below.  It will be proved in a
forthcoming paper that there are only finitely many other Seifert
fibered surgeries on hyperbolic Montesinos knots of length 3.

\section{Preliminaries and examples} \label{sec:2}

Given a submanifold $F$ of a manifold $M$, let $N(F)$ be a regular
neighborhood of $F$ in $M$.  When $F$ has codimension 1 and is
properly embedded, denote by $M|F$ the manifold obtained by cutting
$M$ along $F$.  If $K$ is a knot in $M$, denote by $(M,K,r)$ the
manifold obtained from $M$ by Dehn surgery on $K$ along a slope $r$ on
$\bdd N(K)$.  When $M = S^3$, simply denote $(S^3, K,r)$ by $K(r)$.

A {\em cusped manifold\/} is a compact 3-manifold $M$ with a specified
{\em vertical boundary\/} $\bdd_v M$, which is a disjoint union of
annuli and tori on $\bdd M$.  The closure of $\bdd M - \bdd_v M$ is
the {\em horizontal boundary\/} of $M$, denoted by $\bdd_h M$.  If $M$
is an $I$-bundle over a compact surface $F$ then it has a natural
cusped manifold structure with $\bdd_v M$ the annuli over $\bdd F$.
Conversely, a cusped manifold $M$ is considered an $I$-bundle only if
it is an $I$-bundle with $\bdd_v M$ defined above.  A surface $F$
properly embedded in $M$ with $\bdd F \subset \bdd_h M$ is an {\it
  h-essential surface\/} if it is incompressible, and has no boundary
compressing disk disjoint from $\bdd_v M$.

Let $K$ be a hyperbolic knot in a solid torus $V$ with $\wrap(K)=2$.
Let $D$ be a meridional disk of $V$ intersecting $K$ twice.  Cutting
$V$ along $D$, we obtain a 3-ball $B$.  Let $\tau = K\cap B$ be the
2-string tangle in $B$.  Denote by $X = B - \Int N(\tau)$ the tangle
space.  Clearly $X$ is irreducible, and the hypothesis that $K$ is
hyperbolic implies that $X$ is also atoroidal.  Define the vertical
boundary of $X$ to be $\bdd_v X = \bdd V \cap X$.  Then the horizontal
boundary $\bdd_h X$ is the disjoint union of two copies of once
punctured torus when $\wind(K) = 0$, or a single twice punctured torus
when $\wind(K) = 2$.

Let $Y = N(D\cup K)$ and define $\bdd_v Y = \bdd V \cap Y$.  Then we
can write $V = X \cup _{\eta} Y$, where $\eta: \bdd_h X \to \bdd_h Y$
is a homeomorphism.  The surgery manifold can then be written as
$(V,K,r) = X \cup _{\eta} (Y,K,r)$.

Fix a trivial embedding of $V$ in $S^3$.  Let $K'$ be the core of $S^3
- V$.  If $K$ is a knot in $V$ then $L = K' \cup K$ is a link in
$S^3$.  Conversely, if $L = K' \cup K$ is a two component link in
$S^3$ and $K'$ is trivial then $K$ is a knot in $V = S^3 - \Int
N(K')$.  We use the convention that a trivial circle $K'$ with a dot
represents the component that need to be deleted, so the link $L = K'
\cup K \subset S^3$ represents the pair $(V, K)$ with $V = S^3 - \Int
N(K')$.  The preferred meridian-longitude pair $(m,l)$ of $K$ in $S^3$
(see \cite{Ro}) is then considered the {\it preferred
  meridian-longitude pair\/} of $K$ in $V$.  This sets up a coordinate
system for the slopes on $\bdd N(K)$, so a slope $ql+pm$ is
represented by a rational number $p/q$, or $1/0$ if $(p,q) = (1,0)$.

Let $C$ be the core of $V$.  Fix a meridian-longitude pair $(m_0,l_0)$
of $\bdd V$.  We can re-embed $V$ in $S^3$ by an orientation
preserving homeomorphism $\varphi_n: V \to V$ such that $l_0$ is
mapped to the curve $l_0 + nm_0$ on $\bdd V$.  Denote by $K_n =
\varphi_n(K)$.  Thus $K_n$ is obtained from $K$ by $n$ right hand full
twists along a disk bounded by $K'$.  If $r$ is a slope on $\bdd
N(K)$, let $r_n$ be the corresponding slope $\varphi_n(r)$ on $\bdd
N(K_n)$.  We have $r_n = r + n \times \wind(K)^2$, hence $r_n = r$ if
$\wind(K) = 0$, and $r_n = r + 4n$ if $\wind(K) = 2$.

\bigskip
\leavevmode

\centerline{\epsfbox{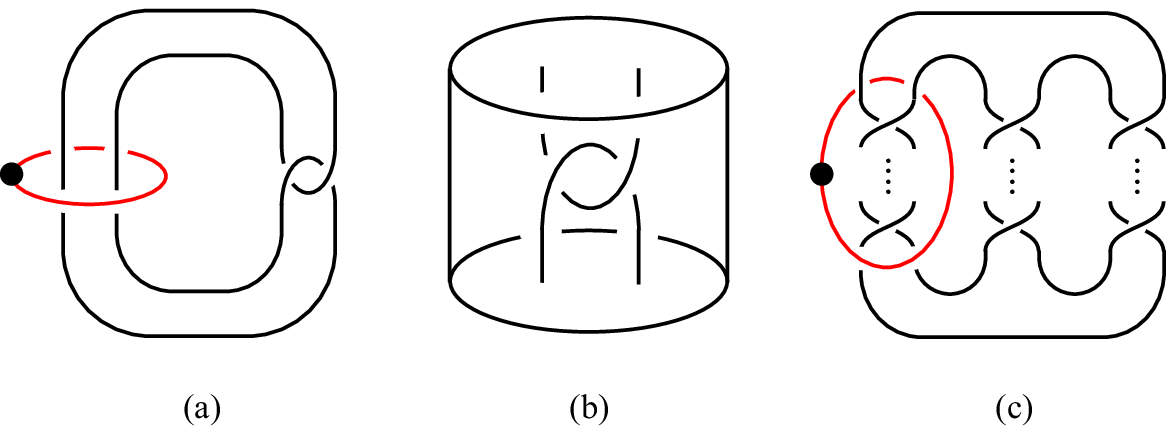}}
\bigskip
\centerline{Figure 2.1}
\bigskip

\begin{example} \label{example:2.1} {\rm (1) Let $K$ be the Whitehead
    knot in $V$ as shown in Figure 2.1(a).  Then $(V,K,r)$ is toroidal
    for $r = 0,4$, and is Seifert fibered for $r=1,2,3$.  See
    \cite[Lemma 7.1]{GW2} and \cite[Lemma 2.3]{BW}.  Cutting $(V,K)$
    along a meridional disk produces a tangle $(B,\tau)$ as shown in
    Figure 2.1(b), which will be called the Whitehead tangle.

    (2) Suppose $K$ has a spanning surface $F$ in $V$ which is a once
    punctured torus or Klein bottle with boundary slope $r$.  Then $F$
    becomes a closed surface $\hat F$ in $(V,K,r)$, which is either a
    Klein bottle or a nonseparating torus.  Since $(V,K,r)$ is
    irreducible \cite{Sch}, the boundary of a regular neighborhood of
    $\hat F$ is incompressible in $(V,K,r)$, hence $(V,K,r)$ is
    toroidal.  In particular, any hyperbolic pretzel knot in solid
    torus as shown in Figure 2.1(c) admits a toroidal surgery along
    the boundary of its pretzel surface. }
\end{example}

\bigskip
\leavevmode

\centerline{\epsfbox{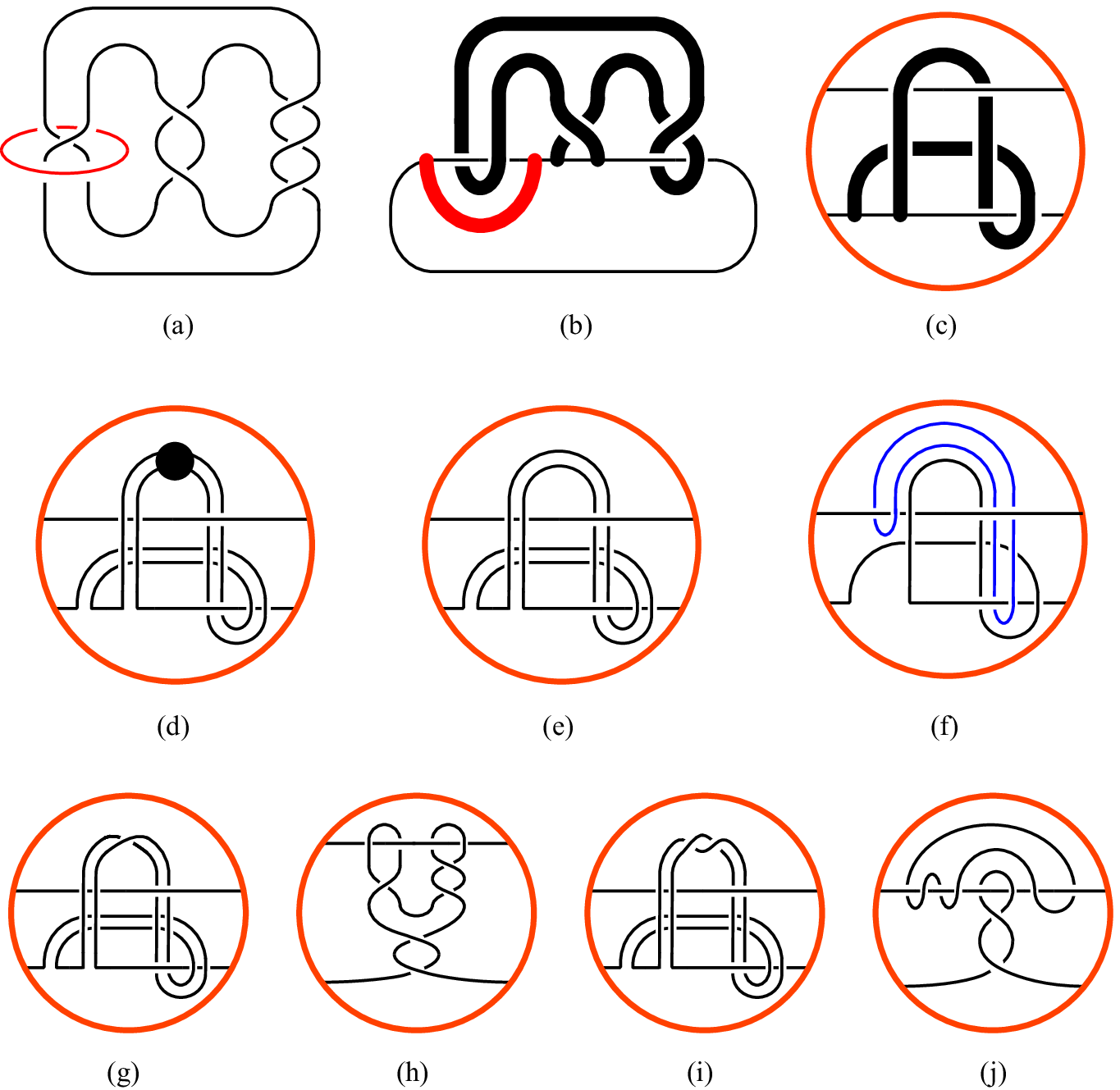}}
\bigskip
\centerline{Figure 2.2}
\bigskip

Note that $K_n(r_n)$ is obtained from $(V,K,r)$ by attaching a solid
torus to $\bdd V$ along the slope $l_0 - n m_0$ on $\bdd V$, hence if
$(V,K,r)$ is hyperbolic then $K_n(r_n)$ is hyperbolic for all but
finitely many $n$.

\begin{prop} \label{prop:2.2} Let $K$ be the knot in solid torus as
  shown in Figure 2.2(a).  Then $(K,V,8)$ and $(K,V,6)$ are toroidal,
  and $(K,V,7)$ is small Seifert fibered with two singular fibers of
  indices $3$ and $5$, respectively.
\end{prop}

\proof Rotating along a horizontal axis of the knot diagram gives
a double branched cover of $(S^3,L)$ over the pair $(S^3, \lambda)$ in
Figure 2.2(b), where $\lambda$ is a pair of arcs represented by the
thick curves in the figure.  The thin circle $C$ in Figure 2.2(b) is the
image of the axis and hence forms the branch set.  The neighborhoods
of the two components of $L = K' \cup K$ in Figure 2.2(a) project to
regular neighborhoods of $\lambda$, which are 3-balls $B_1, B_2$
respectively, where $B_1$ is represented by the lower thick arc in
Figure 2.2(b).  Let $B$ be the closure of the complement of $B_1$, and
let $\tau = C \cap B$.  Then $(B, B_2, \tau)$ can be deformed to that
in Figure 2.2(c) and then further to that in Figure 2.2(d).

Since $V$ is the exterior of $K'$, the above shows that $(V,N(K))$ is
the double branched cover of $(B,B_2)$ branched over $\tau$.  Put
$\tau_2 = \tau \cap B_2$, and denote by $(B, \tau(s))$ the tangle
obtained from $(B, \tau)$ by replacing the subtangle $(B_2, \tau_2)$
with a rational tangle of slope $s$ with respect to certain coordinate
system on $\bdd B_2$, set up so that $\tau_2$ has slope $\infty$.  Let
$r_0$ be the slope on $\bdd N(K)$ which covers a curve of slope $0$ on
$\bdd B_2$.  Then by the Montesinos trick \cite{Mo}, $(V,K,r_0 - s)$ is
then the double branched cover of $B$ branched along the tangle
$\tau(s)$, and $(V,K,r_0-s)$ is Seifert fibered if and only if
$(B,\tau(s))$ is a Montesinos tangle.  

To determine the slope $r_0$, consider the pretzel surface for the
knot $K$ in Figure 2.2(a).  It is a once punctured Klein bottle $F$.
The boundary of $F$ is the pretzel framing $\lambda$, and one can show
that it is a curve of slope $8$ on $\bdd N(K)$ with respect to the
preferred meridian-longitude of $K$.  The projection of $F$ is a disk
$F'$ intersecting the axis at one arc and two individual points, and
the boundary of $F'$ contains an arc $\lambda'$ on $\bdd B_2$ which is
the projection of the above pretzel framing and will be called the {\em
  pretzel framing\/} on $\bdd B_2$.  In Figure 2.2(b) $B_2$ is the
thick dark arc.  Its boundary then contains a pair of arcs connecting
the 4 branch points, called the {\em blackboard framing}.  In our case
these two framings are actually the same because $F'$ has boundary on
the blackboard framing except at the two crossings of the dark curve,
which contributes $-2$ and $2$ respectively to the pretzel framing and
therefore canceled.  One can check that the blackboard framing is
unchanged under the isotopy from Figure 2.2(b) to Figure 2.2(d) and
therefore represents the $0$ slope on $\bdd B_2$.  It follows that the
pretzel slope $8$ is the $r_0$ if we set up the coordinate on $\bdd B$
in the standard way, i.e.\ the horizontal arcs connecting the branch
points represent slope $0$ and the vertical arcs represent $\infty$.
It follows that $(V,K,8-s)$ is the double branch cover of $(B,
\tau(s))$.  In particular, $K(7)$ and $K(6)$ are the double branched
cover of $(B,\tau(1))$ and $(B,\tau(2))$, respectively.

Since $8$ is the pretzel slope, by Example \ref{example:2.1}(2) we see
that $(V,K,8)$ is toroidal.  This can be verified as follows.  The
tangle $\tau(0)$ is shown in Figure 2.2(e), which can be deformed to
that in Figure 2.2(f).  Note that it has a closed component which
bounds a disk $D$ intersecting the other components at two points.
The boundary of a regular neighborhood of $D$ is then a Conway sphere,
which lifts to the incompressible torus in $K(8)$ bounding a twisted
$I$-bundle over the Klein bottle.

The tangle $\tau(1)$ is shown in Figure 2.2(g) and (h).  Without fixing
the endpoints of the strings on the outside sphere it is equivalent to
the $(-1/3,\, -1/5)$ Montesinos tangle.  Hence its double branched
cover $(V,K,7)$ is a small Seifert fiber space with two singular fibers
of indices $3$ and $5$, respectively.

The tangle $\tau(2)$ is shown in Figure 2.2(i), which deforms to that
in Figure 2.2(j).  There is an obvious Conway sphere bounding a
$(1/2,\, -1/4)$ Montesinos tangle, and its outside is not a product,
hence it lifts to an essential torus in $(V,K,6)$, bounding a small
Seifert fiber space with two singular fibers of indices $2$ and $4$
respectively.  
\qed \bigskip

The following result shows that each pretzel knot of type
$(-2,3,2n+1)$ admits at least two Seifert fibered surgeries, with
slopes $6+4n$ and $7+4n$.  In particular, when $n = 3$ it gives the
well known results that $18$ and $19$ surgeries on the $(-2,3,7)$
pretzel knot are lens spaces [FS].  Denote by $M(r_1, r_2, r_3)$ the
closed 3-manifold which is the double branched cover of $S^3$ with
branch set the Montesinos link $K(r_1, r_2, r_3)$.  To make the
statement simple, we do allow $r_3 = 0$ in this theorem, in which case
$K(r_1, r_2, r_3)$ is actually the connected sum of two 2-bridge
knots, and $M(r_1, r_2, r_3)$ is reducible.

\begin{cor} \label{cor:2.3} Let $K_n$ be the $(-2,\,3,\, 2n+1)$
  pretzel knot in $S^3$.  Then $K_n(7+4n) = M(-1/3,\, 3/5,\,
  1/(n-2))$, and $K_n(6+4n) = M(1/2,\, -1/4,\, 2/(2n-5))$.  In
  particular, they are Seifert fibered manifolds for all $n$, except
  that when $n=2$, $K_2(15) = M(-1/3,\, 3/5,\, 1/0)$ is reducible.
\end{cor}

\proof Let $r_n = r + 4n$, where $r = 6, 7$.  Recall that $K_n$ is
obtained from $K \subset V \subset S^3$ by $n$ right hand full twists
along a meridian of $V$, so $K_n(r_n)$ is obtained from $(V,K,r)$ by
attaching a solid torus $V'$ on the outside so that a meridian of $V'$
is attached to the curve $\lambda = l_0 - n m_0$.  By the Montesinos
trick, $K_n(r_n)$ is the double branched cover of $S^3$ along the link
$L$ obtained from $(B,\tau(8-r))$ by attaching a rational tangle $(B',
\tau')$ to the outside of $B$.

To calculate the slope of $(B',\tau')$, note that $m_0$ and $l_0$
projects to curves of slope $0/1$ and $1/0$, respectively.  One can
then check that the curve $\lambda$ projects to a curve $\lambda'$ of
slope $-1/n$ on $\bdd B$.  Since the map $\bdd B' \to \bdd B$ is
orientation reversing, $\lambda'$ is of slope $1/n$ on $\bdd B'$.  We
may assume that $\lambda$ has been isotoped to bound a meridian disk
$D$ in $V'$ which is disjoint from the branch axis.  Then $\lambda'$
bounds a disk in $B'$ disjoint from the tangle strings.  It follows
that $(B', \tau')$ is of slope $1/n$.

For $r=7$, the tangle $(B, \tau(1))$ in Figure 2.2(h) is a Montesinos
tangle of length 2, and $K_n(7+4n)$ is the double branched cover of
the link obtained by attaching a $1/n$ tangle to the outside of $(B,
\tau(1))$, which one can check is the link $K(-1/3,\, 3/5,\,
1/(n-2))$.  Hence $K_n(7+4n) = M(-1/3,\, 3/5,\, 1/(n-2))$.  It is a
Seifert fiber space (possibly a lens space) unless $n=2$, which gives
the reducible 15 surgery on the $(3,5)$ torus knot.  For $r=6$, we
note that the union of $(B', \tau')$ and the tangle $(B, \tau(2))$ in
Figure 2.2(j) is the Montesinos link $K(1/2,\, -1/4,\, 2/(2n-5))$,
hence $K_n(6+4n) = M(1/2,\, -1/4,\, 2/(2n-5))$, which is a small
Seifert fiber space for any $n$.  \qed

\section{Surgery on $K \subset V$ with $\wrap(K)=2$}  \label{sec:3}

Throughout this section we will assume that $K\subset V$ is a
hyperbolic knot with $\wrap (K) = 2$, intersecting a meridian disk
$D$ of $V$ twice.  Recall that $Y = N(D \cup K)$, $(B, \tau)$ is the
tangle obtained by cutting $(V,K)$ along $D$, and $X = V - \Int (Y)
= B - \Int N(\tau)$.  Let $r$ be a nontrivial slope such that
$(V,K,r)$ is nonhyperbolic.  Denote by $K_r$ the dual knot in
$(V,K,r)$ and $(Y,K,r)$.

\begin{lemma} \label{lemma:3.1}
Suppose $\wrap(K)=2$ and $\bdd_h X$ is incompressible in
  $X$.  If $X$ is an $I$-bundle with $\bdd_h X$ the $\bdd I$-bundle
  then $K$ is the Whitehead knot in $V$.
\end{lemma}

\proof If $\wind(K) = 0$ then $\bdd_h X$ is the disjoint union of two
copies of once punctured torus.  Hence the hypothesis above implies
that $X$ is a product $Q \times I$, where $Q$ is a once punctured
torus, and $\bdd_v X = \bdd Q \times I$.  Recall that $X = B - \Int
N(\tau)$.  Let $\tau_1, \tau_2$ be the two strands of $\tau$.  Adding
$N(\tau_1)$ to $X$ produces a $D \times I$ with a 1-handle $H_1$
attached to $D\times 1$, and $\tau_1$ is the core of $H_1$.
Similarly $N(\tau_2)$ can be considered as a 2-handle attached to the
solid torus $X\cup N(\tau_1)$.  Since the result is a 3-ball, the core
of the 2-handle $N(\tau_2)$ intersects the meridian of $X\cup
N(\tau_1)$ at a single point.  It is now clear that $\tau = \tau_1
\cup \tau_2$ is the tangle shown in Figure 2.1(b), hence $K$ is a
Whitehead knot in $V$.

If $\wind(K) = 2$ then $\bdd_h X$ is a twice punctured torus, hence if
$X$ is an $I$-bundle then it must be a twisted $I$-bundle over a once
punctured Klein bottle $P$, so we can properly embed $P$ in $X \subset
B \subset S^3$.  This is impossible because the union of $P$ and a
disk on $\bdd B$ would be a closed Klein bottle embedded in $B^3$.
\qed
\bigskip

An isotopy class $[\alpha]$ of a nontrivial simple closed curve
$\alpha$ on $\bdd_h X$ is called an {\it annular slope\/} if $\alpha$
is not parallel to a boundary component on the surface $\bdd_h X$, and
there is an h-essential annulus $A$ in $X$ with $\alpha$ as a boundary
component.  Note that it is possible that the other boundary component
of $A$ could be a boundary parallel curve on $\bdd_h X$ and hence
would not be an annular slope.

\begin{lemma} \label{lemma:3.2}
Suppose $\bdd_h X$ is incompressible and $X$ is
  not an $I$-bundle.  Then there is a non-separating curve $\alpha$ on
  each component of $\bdd_h X$ which is disjoint from any annular
  slope of $\bdd_h X$ up to isotopy.
\end{lemma}

\proof Let $(W, \bdd_h W)$ be the characteristic pair of the pair $(X,
\bdd_h X)$, as defined in \cite{JS}.  Then $\bdd_h W = W \cap \bdd_h
X$ is a subsurface of $\bdd_h X$, and each boundary component of
$\bdd_h W$ is a nontrivial curve on $\bdd_h X$.  By the definition of
characteristic pair, $\bdd_h W$ contains all annular slopes on $\bdd_h
X$ up to isotopy.

First assume $\wind(K)=0$, so each component $F$ of $\bdd_h X$ is a
once punctured torus.  It is easy to see that if the result is false
then some component $G$ of $\bdd_h W \cap F$ is {\it full\/} in the
sense that $F - G$ is in a collar of $\bdd F$; hence it is a once
punctured torus.  Let $W_0$ be the component of $W$ containing $G$.
Since $G$ is not a double cover of any other surface, $W_0$ must be a
trivial $I$-bundle, so $\bdd_h W_0 - G$ is also a once punctured
torus, which must be on $\bdd_h X - F$.  By Lemma \ref{lemma:3.1} $X$
is not an $I$-bundle, hence $A = \bdd_v W_0$ is an essential annulus
in $X$, cutting off a compact 3-manifold $M$ with $\bdd M$ a single
torus.  Since $X$ is atoroidal and $\bdd_h X$ is incompressible, we
see that $\bdd M$ must be compressible inside of $M$, so $X$ being
irreducible (since $K$ is hyperbolic) implies that $M$ is a solid
torus.  Since $X$ is not an $I$-bundle, we see that $A$ runs at least
twice along the longitude of $M$.  It follows that the union of $A$
and an annulus in $Y$ parallel to $\bdd_v Y$ forms an essential torus
in $V-\Int N(K)$, contradicting the assumption that $K$ is a
hyperbolic knot in $V$.

Now assume $\wind(K)=2$.  Then $\bdd_h X$ is a twice punctured torus,
so if $\bdd_h X - \bdd_h W$ does not contain a nonseparating curve
then some component $G$ of $\bdd_h W$ is a once or twice punctured
torus.  Let $W_0$ be the component of $W$ containing $G$.  Since
$\bdd_h X$ has genus one, there is no room for another copy of $G$,
hence $W_0$ must be a twisted $I$-bundle over a once punctured Klein
bottle $R$.  In particular, $\bdd_h W_0$ must be a twice punctured
torus, so $\bdd_h X - \bdd_h W_0$ is a pair of annuli.  We can then
extend an embedding of $R$ in $W_0$ to an embedding of $R$ in $X$ with
$\bdd R \subset \bdd_v X \subset \bdd B$.  The union of $R$ with a
disk on the boundary of $B$ would then be a closed Klein bottle
embedded in the 3-ball $B$, which is impossible.  \qed

\begin{lemma} \label{lemma:3.3}
Suppose $K\subset V$ is a hyperbolic knot with $\wrap(K)
  = 2$.  

  (1) $\bdd_h Y$ is incompressible in $(Y,K,r)$ for all nontrivial
  $r$.

  (2) If $r$ is an integral slope then $(Y,K,r)$ is an $I$-bundle with
  $\bdd_v Y$ as its vertical surface.

  (3) If $r$ is a nontrivial non-integral slope then any h-essential
  annulus $Q$ in $(Y,K,r)$ can be isotoped to be disjoint from the
  dual knot $K_r$.
\end{lemma}

\proof Recall that $Y = N(D \cup K)$, where $D$ is a meridian disk
of $V$ intersecting $K$ twice.  Let $D_1$ be a meridian disk of $K$ in
$Y$ with $\bdd D_1 \subset \bdd_h Y$, and let $Y_1 = N_1(D_1 \cup K)$
be a smaller regular neighborhood of $D_1 \cup K$ such that $Y_1 \cap
\bdd Y = \bdd D_1 \times I$.  Then the frontier of $Y_1$ is an annulus
$A$, cutting $Y$ into $Y_1$ and another manifold $W$.  When
$\wind(K)=0$ $W$ is a product $T_1 \times I$, where $T_1$ is a once
punctured torus; when $\wind(K)=2$ the manifold $W$ is a twisted
$I$-bundle over a Klein bottle.  In either case $W$ is an $I$-bundle
with $\bdd_v Y$ as its vertical boundary.  Note that $\bdd_h W$ is
incompressible, and $A$ is an annulus on $\bdd_h W$, which is
incompressible in $W$.

It is clear that $Y_1$ is a solid torus with $K$ as a core, hence $V'
= (Y_1, K, r)$ is a solid torus for all $r$.  When $r$ is an integral
slope $A$ runs along the longitude of $V'$ once, hence $(Y,K,r) = W
\cup_A (Y_1,K,r)$ is homeomorphic to the $I$-bundle $W$ with $\bdd_v
Y$ preserved.  When $r$ is a nontrivial non-integral slope $A$ runs
along the longitude of $V'$ more than once.  By a standard innermost
circle outermost arc argument one can show that $\bdd_h Y$ is
incompressible in $(Y,K,r)$.

If $Q$ is an h-essential annulus in $(Y,K,r)$ then it can be isotoped
so that $Q \cap A$ has no arc component, so $Q \cap (Y_1, K, v)$ is a
set of incompressible annuli, which can then be isotoped to be
disjoint from $K_r$.  \qed

\begin{lemma} \label{lemma:3.4}
Suppose $K$ is a hyperbolic knot in $V$ with
  $\wrap(K)=2$, $K$ is not the Whitehead knot, and $\bdd_h X$ is
  incompressible in $X$.

(1) $(V, K, r)$ is irreducible, $\bdd$-irreducible, and is not Seifert
    fibered.

(2) If $r$ is not in integral slope then $(V,K,r)$ is hyperbolic.
\end{lemma}

\proof (1) The irreducibility and $\bdd$-irreducibility follows from
\cite{Be, Ga2, Sch}.  It also follows from Lemma \ref{lemma:3.3}
because $\bdd_h X$ is an essential surface in $(V,K,r)$ and there is
no reducing sphere or compressing disk of $\bdd V$ disjoint from
$\bdd_h X$.

Suppose $(V,K,r)$ is Seifert fibered.  By \cite{Wa} any incompressible
surface in a Seifert fibered space is either horizontal or vertical.
Since $\bdd_h X$ is not an annulus or torus, it must be horizontal,
so both $X$ and $Y(r)$ are $I$-bundles with $\bdd_h X = \bdd_h Y$ as
their horizontal surface.  On the other hand, by Lemma \ref{lemma:3.1}
this is impossible unless $K$ is the Whitehead knot in $V$, which has
been excluded.

(2) If $T$ is an essential torus in $(V,K,r)$ then it must intersect
$\bdd_h X$ because both $(Y,K,r)$ and $X$ are atoroidal.  Using a
standard cut and paste argument one can show that $T$ can be isotoped
so that each component of $T \cap X$ and $T \cap (Y,K,r)$ is an
h-essential annulus.  If $r$ is not an integral slope then by Lemma
\ref{lemma:3.3}(3) the annuli $T \cap (Y,K,r)$ can be isotoped to
be disjoint from $K_r$, so $T$ would be an essential torus in $V-K$,
contradicting the assumption that $K$ is hyperbolic.  \qed
\bigskip

A curve $\alpha$ on a surface $F$ is {\it orientation preserving\/} if
the orientation of $F$ does not change when traveling through
$\alpha$.  Alternatively, $\alpha$ is orientation preserving if its
regular neighborhood is an annulus, not a M\"obius band.

\begin{lemma} \label{lemma:3.5}
  Up to isotopy there are exactly two orientation preserving essential
  simple closed curves on a Klein bottle $F$.
\end{lemma}

\proof Let $\a$ be a curve cutting $F$ into an annulus $A$.  Suppose
$\b$ is another orientation preserving essential simple closed curve,
which intersects $\alpha$ minimally but is not isotopic to $\alpha$.
Then $\a$ cuts $\b$ into a set of essential arcs $C$ on $A$.  One can
show that $C$ has exactly two components as otherwise $\beta$ would
either be orientation reversing (if $|\b \cap \a|$ is odd) or contains
more than one components (if $|\b \cap \a| > 2$).  Therefore any other
such curve is obtained from $\b$ by Dehn twists along $\a$.  It is
easy to check that Dehn twisting $\b$ once along $\a$ produces a curve
isotopic to $\b$, hence the result follows. 
\qed 
\bigskip

We note that the two curves $\a, \b$ in the above proof are
essentially different as $\a$ cuts $F$ into an annulus while $\b$ cuts
$F$ into two M\"obius bands.  Define an orientation preserving essential
simple closed curve $\gamma$ on a surface $F$ to be {\it of type I or
  type II\/} according to whether $F|\gamma$ is orientable or not.
Thus the curve $\alpha$ above is of type I and $\beta$ of type II.  If
$C$ is an annular slope of a twisted $I$-bundle $W$ over $F$ and
$\varphi: W \to F$ the $I$-fibration, then up to isotopy $C$ is a
boundary component of a vertical annulus.  Hence we can define $C$ to
be {\it of type I or type II\/} according to whether $\varphi(C)$ is
of type I or II on $F$.

Given a compact surface $F$, denote by $\hat F$ the closed surface
obtained from $F$ by capping off each boundary component with a disk.
Two curves $C_1, C_2$ on $F$ are {\em weakly equivalent}, denoted by
$C_1 \sim C_2$, if they are isotopic on $\hat F$.

Let $W$ be a twisted $I$-bundle over a once punctured Klein bottle
$R$.  Let $\hat W$ be the manifold obtained by attaching a 2-handle to
$W$ along $\bdd_v W$.  Then $\hat W$ is a twisted $I$-bundle over the
Klein bottle $\hat R$.  Denote by $\varphi: \bdd_h W \to \bdd \hat W$
the inclusion map.  Let $C_1$ be a type I annular slope on $\bdd_h W$.
There are infinitely many annular slopes on $F$ that intersect $C_1$
essentially, but the following shows that these are all weakly
equivalent.

Denote by $I(C_1, C_2)$ the algebraic intersection number between two
curves $C_1, C_2$, which is well defined up to sign on any orientable
surfaces.

\begin{lemma}  \label{lemma:3.6}
  Let $C_1$ be a type I annular slope on $\bdd_h W$.  Then there is
  a type II annular slope $C_2$ on $\bdd_h W$ intersecting $C_1$ at
  a single point, such that if $C_3$ is an annular slope on $\bdd_h W$
  then it is either weakly equivalent to $C_2$, or isotopic on
  $\bdd_h W$ to a curve disjoint from $C_1$.  In particular, $I(C_3,
  C_i) = 0$ for some $i=1,2$.
\end{lemma}

\proof Fix $I$-bundle structures of $W$ and $\hat W$ and let $\rho: W
\to R$ and $\hat \rho: \hat W \to \hat R$ be the projection maps.  We
may assume that $C_1$ is a boundary component of a vertical annulus
$A_1$.  Then $\rho(A_1) = \rho(C_1) = \alpha_1$ is a type I curve on
$R$.  Let $\a_2$ be a type II curve on $R$ intersecting $\alpha_1$
minimally at two points as given in the proof of Lemma
\ref{lemma:3.5}, and let $C_2$ be a boundary component of
$\rho^{-1}(\a_2)$.  The two intersection points of $\a_1 \cap \a_2$
lifts to two points of $C_1 \cap \rho^{-1}(\a_2)$, one on each
component of $\rho^{-1}(\a_2)$.  Hence $C_1$ intersects $C_2$ at a
single point.

Assume $C_3$ is another annular slope on $F$.  We may assume it is a
boundary component of a vertical annulus $A_3$, so $\alpha_3 =
\rho(A_3)$ is an orientation preserving essential simple closed curve
on $R$.  By Lemma \ref{lemma:3.5} $\a_3$ is isotopic to either $\a_1$
or $\a_2$ on $\hat R$.  Cutting $R$ along $\a_1$ produces a once
punctured annulus, hence it is easy to see that if $\a_3$ is isotopic
to $\a_1$ on $\hat R$ then it is also isotopic to $\a_1$ on $R$, in
which case $C_3$ can be isotoped to be disjoint from $C_1$.  If $\a_3$
is isotopic to $\a_2$ then $C_3$ is isotopic to a component of $\hat
\rho^{-1}(\a_2)$ on $\bdd \hat W$.  Since the two components of $\hat
\rho^{-1}(\a_2)$ are parallel to each other, we see that $C_3$ is
isotopic to $C_2$ on $\bdd \hat W$ and hence is weakly equivalent to
$C_2$ on $\bdd W$.  \qed

\begin{lemma} \label{lemma:3.7} Suppose $\wind(K) = 2$.  Let $r$ be an
  integral slope and $K_r$ the dual knot in $(Y,K,r)$.  Let $\a$ be a
  simple closed curve on $\bdd_h Y$ which is isotopic to $K_r$ in
  $(Y,K,r)$, let $\b$ be an annular slope on $\bdd_h Y$ intersecting
  $\a$ essentially at one point, as given in Lemma \ref{lemma:3.6}.
  Suppose $s \neq r$ is another integral slope on $\bdd N(K)$.  Then
  $\beta$ is not weakly equivalent to an annular slope of $\bdd_h Y$
  in $(Y,K,s)$.
\end{lemma}

\proof By the proof of Lemma \ref{lemma:3.3}, $(Y,K,r)$ is obtained
from an $I$-bundle $W$ over $R$ by attaching a solid torus $V'$ along
a longitudinal annulus of $V'$, and $K_r$ is the core of $V'$.  It is
easy to see that $\a$ is a type $I$ annular slope.  The identity map
of $W$ extends to a homeomorphism $\psi: (Y,K,s) \to (Y,K,r)$, and the
restriction of $\psi$ on $\bdd Y$ is a Dehn twist $\tau_{\alpha}^n$
along $\a$, where $n = s-r \neq 0$.  In particular, the curve $\beta$
is mapped to $\beta' = \tau_{\alpha}^n (\beta)$ on $\bdd_h (Y,K,r)$.
We have $|I(\b', \a)| = 1$ and $|I(\b', \b)| = |n| \neq 0$; hence by
Lemma \ref{lemma:3.6} $\b'$ is not weakly equivalent to an annular
slope of $\bdd_h (Y,K,r)$.  Since the homeomorphism $\psi: (Y,K,s) \to
(Y,K,r)$ maps $\beta$ to $\beta'$, it follows that $\beta$ is not
weakly equivalent to an annular slope of $\bdd_h (Y,K,s)$.  \qed

\begin{thm} \label{thm:3.8} 
  Suppose $K$ is a hyperbolic knot with $\wrap(K) = 2$ in a solid
  torus $V$, $K$ is not the Whitehead knot, and $\bdd_h X$ is
  incompressible in $X$.  Then $K$ admits at most one exceptional
  surgery $(V,K,r)$, which must be a toroidal surgery and $r$ an
  integral slope.
\end{thm}

\proof By Lemma \ref{lemma:3.4} $(V,K,r)$ is irreducible and not a
solid torus or small Seifert fiber space, and it is also atoroidal
when $r$ is not an integral slope.  Hence we need only show that $K$
admits at most one integral toroidal surgery.

First assume $\wind(K) = 2$.  By Lemma \ref{lemma:3.2} there is a
nonseparating curve $\gamma$ on $\bdd_h X$ which is disjoint from all
annular slopes of $X$ up to isotopy.  Suppose $(V,K,r)$ is toroidal.
Let $\alpha, \beta$ be the annular slopes on $\bdd_h Y = \bdd_h
(Y,K,r)$ as given in Lemma \ref{lemma:3.7}.  Let $q: \bdd_h X \to
\bdd_h Y$ be the gluing map.  We claim that $q(\gamma) \sim \beta$.

Let $T$ be an essential torus in $(Y,K,r)$ intersecting $\bdd X$
minimally.  Since $\bdd_h X$ is incompressible, each component of $A_1
= T \cap X$ and $A_2 = T\cap (Y,K,r)$ is an h-essential annulus.  In
particular, each boundary component of $A_1$ is either an annular
slope of $X$ or boundary parallel, hence by the above we may assume
$\gamma \cap \bdd A_1 = \emptyset$, so $q(\gamma) \cap \bdd A_2 =
\emptyset$.  On the other hand, each boundary component of $A_2$ is
either parallel to a boundary component of $\bdd_h Y$, or is an
annular slope of $(Y,K,r)$; hence by Lemma \ref{lemma:3.6} we may
assume that it is either disjoint from $\a$ or weakly equivalent to
$\b$.  If no component of $\bdd A_2$ is weakly equivalent to $\b$ then
$\bdd A_2$ can be isotoped to be disjoint from $\a$, hence $T$ can be
isotoped to be disjoint from $K_r$ because by definition $K_r$ is
isotopic to $\a$.  This contradicts the assumption that $K$ is
hyperbolic.  Therefore at least one component $C$ of $\bdd A_2$
satisfies $C \sim \b$; in particular, it is nonseparating.  Since $C$
and $q(\gamma)$ are disjoint and they are both nonseparating curves on
the punctured torus $\bdd_h Y$, we have $C \sim q(\gamma)$, hence the
claim $q(\gamma) \sim \b$ follows.

Now if $s$ is another toroidal slope of $K$ then by the above, $\beta
\sim q(\gamma)$ is also an annular slope on $\bdd_h (Y,K,s)$, which
contradicts Lemma \ref{lemma:3.7}, completing the proof for the case
of $\wind(K) = 2$.

The proof for the case of $\wind(K)=0$ is similar.  In this case
$(Y,K,r)$ is a product $F \times I$.  Let $F_i = F \times i$ for
$i=0,1$. Let $q: \bdd_h X \to \bdd_h Y = F_0 \cup F_1$ be the gluing
map, and $G_i = q^{-1}(F_i)$.  By Lemma \ref{lemma:3.2} there is a
nonseparating curve $\gamma_i$ on $G_i$ which is disjoint from all
annular slopes of $X$ up to isotopy.  Let $\gamma'_i$ be the curve
$q(\gamma_i)$ on $G_i$.  We claim that $\gamma'_0, \gamma'_1$ cobound
an annulus and hence is homologous in $(Y,K,r)$.  

Let $T$ be an essential torus in $(V,K,r)$.  As above, the
hyperbolicity of $K$ implies that there is a component $A'$ of $A_1 =
T \cap (Y,K,r)$ which cannot be isotoped off $K_r$.  Let $\beta_i =
A'\cap F_i$.  Since each side of $A'$ must be adjacent to an essential
annulus in $X$, we see that $q^{-1}(\beta_i)$ is an annulus slope on
$G_i$.  Since $G_i$ is a once punctured torus, any annular slope on
$G_i$ is a nonseparating curve disjoint from $\gamma_i$ and therefore
must be isotopic to $\gamma_i$.  If follows that $\gamma'_i$ is
isotopic to $\beta_i$ on $F_i$, hence $A'$ can be isotoped to have
$\bdd A' = \gamma'_0 \cup \gamma'_1$, and the claim follows.

For the same reason, if $s$ is another integral toroidal slope of $K$
then $\gamma'_0$ and $\gamma'_1$ must also be homologous in $(Y,K,s)$.
On the other hand, as in the proof of Lemma \ref{lemma:3.7}, there is a
homeomorphism $\psi: (Y,K,r) \to (Y,K,s)$ which is the identity map on
$F_0$ and the Dehn twist map $\tau_{\alpha}^{n}$ on $F_1$, where $n =
r-s \neq 0$ and $\alpha$ is the curve on $F_1$ isotopic to $K_r$ in
$(Y,K,r)$.  Since $A'$ cannot be isotoped off $K_r$, $\gamma'_1$ has
essential intersection with $\alpha$, hence $\psi(\gamma'_1)$ is not
homologous to $\gamma'_0$ in $(Y,K,s)$ if $s \neq r$, a contradiction.
\qed

\section{Surgery on $K_n$}  \label{sec:4}

As in Section 1, define $K_n = \varphi_n(K)$ and $r_n = \varphi_n(r)$,
where $\varphi_n: V \to S^3$ is the composition of the standard
embedding of $V$ into $S^3$ with $n$ full right hand twists along a
meridian disk of $V$, and $r$ is a slope of $K$.  If $(V,K,r)$ is a
small Seifert fiber manifold then $K_n(r_n)$ is either small Seifert
fibered or reducible, hence is always nonhyperbolic.  In general, if
$(V,K,r)$ is toroidal then it is possible that $K_n(r_n)$ may be
hyperbolic for infinitely many $n$; however, this will not happen if
$\wind(K) = 2$.  The main result of this section shows that if
$\wrap(K) = 2$ and $(V,K,r)$ is toroidal then either $K_n(r_n)$ is
toroidal for all but at most three $n$, or it is atoroidal and
nonhyperbolic for all $n$.  In particular, it can be hyperbolic for at
most three $n$.  See Theorem \ref{thm:4.3} for more details.

Let $D$ be a meridional disk of $V$ with $n_1 = |D \cap K| = 2$, and
$T$ an essential torus in $(V,K,r)$ such that $n = n_2 = |T \cap K_r|$
is minimal.  Let $E(K)$ be the knot exterior $V - \Int N(K)$.  Denote
by $Q_1$ the punctured disk $D \cap E(K)$, and by $Q_2$ the punctured
torus $T \cap E(K)$.  Considering the
disks $D \cap N(K)$ and $T \cap N(K_r)$ as fat vertices, and the arc
components of $Q_1 \cap Q_2$ as edges, we obtain graphs $\Gamma_1,
\Gamma_2$ on $D$ and $T$, respectively, with $n_i$ vertices on
$\Gamma_i$.  Denote by $m$ the meridional slope of $K$, and by $\Delta
= \Delta(m, r)$ the distance (i.e.\ the geometric intersection number)
between $m$ and $r$.  By \cite{GL2} we have $\Delta \leq 2$.  Each
boundary component of $Q_1$ intersects each component of $Q_2$ at
$\Delta$ point; hence each vertex of $\Gamma_1$ has valence $n
\Delta$, and each vertex of $\Gamma_2$ has valence $n_1 \Delta = 2
\Delta$.

The above are standard set up for intersection graphs of exceptional
Dehn surgeries.  We refer the readers to \cite{CGLS, GW1} for standard
terms and basic results related to intersection graphs, such as
Scharlemann cycles, extended Scharlemann cycles, and signs of
vertices.  In particular, the minimality of $n$ and $\wrap(K)=2$ imply
that there is no trivial loops in $\Gamma_i$, so $\Gamma_1$ is a set
of $n \Delta$ parallel edges.  Each vertex of $\Gamma_i$ has a sign.
An edge of $\Gamma_i$ is a positive edge if the two vertices on its
endpoints are of the same sign.  There is a one to one correspondence
between edges of $\Gamma_1$ and $\Gamma_2$.  The Parity Rule of
\cite[P279]{CGLS} says that an edge is positive on one graph if any
only if it is negative on the other.  If $\wind(K) = 2$ then both
vertices of $\Gamma_1$ are positive, hence all edges on $\Gamma_1$ are
positive edges, and all edges on $\Gamma_2$ are negative edges.
Similarly if $\wind(K) = 0$ then all edges of $\Gamma_1$ are negative
and all edges of $\Gamma_2$ are positive.

\begin{lemma} \label{lemma:4.1}
Suppose $K\subset V$ has $\wind(K) = \wrap(K) = 2$.  If
  $(V,K,r)$ is toroidal then it contains an essential torus $T$ such
  that

(1) If $n > 2$ then $\Gamma_1$ has no extended Scharlemann cycle;

(2) $n = 2$ or $4$; 

(3) $\Delta = 1$; 

(4) $T$ bounds a small Seifert fiber space.
\end{lemma}

\proof (1) This is \cite[Lemma 2.9]{BZ} or \cite[Theorem 3.2]{GL1}.
An extended Scharlemann cycle can be used to find an essential torus
in $(V,K,r)$ which has fewer intersection with $K_r$, contradicting
the minimality of $n$.

(2) The parity rule implies that $n$ must be even as otherwise there
would be an edge in $\Gamma_1$ with the same label on its two
endpoints, so it would be a positive edge on both graphs.  If $n > 4$
then the $n$ parallel positive edges of $\Gamma_1$ contain an extended
Scharlemann cycle, contradicting (1).  See \cite[Lemma 1.4]{Wu98}.

(3) Assume $\Delta \geq 2$.  If $n = 4$ then $\Gamma_1$ has an
extended Scharlemann cycle, a contradiction.  Hence we may assume $n =
2$.  By \cite[Lemma 2.1]{Go}, no two edges can be parallel on both
graphs, so we must have $\Delta = 2$, and the four edges of $\Gamma_2$
are mutually non-parallel on $\Gamma_2$.  A disk face $E$ of
$\Gamma_2$ then has four edges.  Now cut $V$ along $D$, let $D_1, D_2$
be the two copies of $D$ on $B = V|D$, and let $\tau = \tau_1 \cup
\tau_2 = K \cap B$ be the two strings of $K$ in $B$.  Then the
neighborhood of $D_1 \cup D_2 \cup \tau_1 \cup \tau_2$ is a solid
torus $V'$ in $B$.  The boundary curve of $E$ runs four times along
$\tau$, twice along each $\tau_i$.  Since all segments of $\bdd E$ on
$D$ and $\bdd N(\tau)$ are essential arcs, we see that $\bdd E$
intersects a meridian of $\tau_i$ twice in the same direction, hence
$V'\cup N(E))$ is a punctured projective space in the 3-ball $B$,
which is impossible.

(4) We now have $\Delta=1$ and $n=2$ or $4$.  The edges of $\Gamma_1$
form one or two Scharlemann cycles, according to whether $n=2$ or $4$.
See Figure 4.1.  By \cite[Lemma 2.5.2]{CGLS} the essential torus $T$ is
separating in $(V,K,r)$.  It cuts $(V,K,r)$ into two regions; the one
containing $\bdd V$ is called the {\em white region}, and the other
one the {\em green region}.  From Figure 4.1 we can see that the
Scharlemann disk $G$ bounded by a Scharlemann cycle $e_1\cup e_2$ is
in the green region since there is no extended Scharlemann cycle.
When shrinking each fat vertex of $\gamma_2$ to a point, $e_1 \cup
e_2$ becomes a loop on $T$, which must be essential by \cite[Lemma
2.8]{BZ}.  Let $H$ be the part of $N(K_r)$ in the green region.  Then
$N(T \cup H \cup G)$ has two torus boundary component, and the one
$T'$ inside the green region has fewer intersection with $K_r$.  By
the choice of $T$, this $T'$ must bound a solid torus $V'$.  The green
region is now the union of two solid tori $V'$ and $V'' = N(H\cup E)$
with $V' \cap V''$ an annulus, hence $T$ being essential implies that
$V' \cup V''$ is a small Seifert fiber space bounded by $T$, whose
orbifold is a disk with two cone points.  \qed

\bigskip
\leavevmode

\centerline{\epsfbox{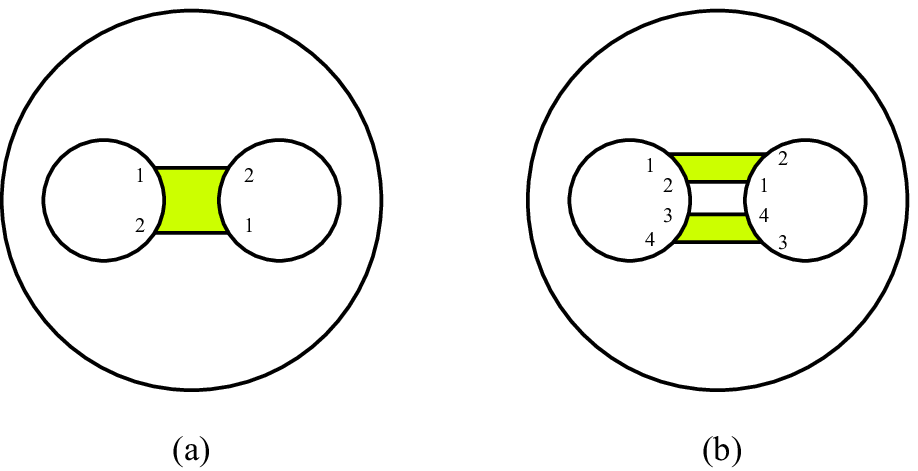}}
\bigskip
\centerline{Figure 4.1}
\bigskip

\begin{lemma} \label{lemma:4.2}
Suppose $K\subset V$ has $\wrap(K) = 2$.  If $(V,K,r)$
  is toroidal then $\Delta = 1$, and it contains an essential torus
  $T$ such that either $T$ is nonseparating or it bounds a small
  Seifert fiber space.
\end{lemma}

\proof The conclusion holds if $(V,K,r)$ contains a nonseparating
torus, so we may assume that all essential tori in $(V,K,r)$ are
separating.  If $(V,K,r)$ contains a Klein bottle $F$ then $T = \bdd
N(F)$ must be incompressible as otherwise $(V,K,r)$ would be
reducible.  $T$ is also not parallel to $\bdd V$, otherwise $(V,K,r) =
N(F)$ would be atoroidal.  Hence $T$ is an essential torus bounding
the small Seifert fiber space $N(F)$ and the result follows.
Therefore we may also assume that $(V,K,r)$ contains no Klein bottle.

The case $\wind(K) =2$ is covered by Lemma \ref{lemma:4.1}, so we
assume $\wind(K) = 0$.  The two vertices of $\Gamma_1$ on $D$ are now
antiparallel, so all edges of $\Gamma_1$ are negative.  By \cite[Lemma
2.3(1)]{GW3}, if $\Gamma_1$ has more than $n$ parallel negative edges
then $T$ would be nonseparating, contradicting the assumption above.
Hence we must have $\Delta = 1$.

By the proof of \cite[Lemma 2.2(3)]{GW3}, the $n$ edges form mutually
disjoint essential cycles of equal length on $\Gamma_2$.  All vertices
on the same cycle are parallel; since $T$ is separating, the number of
positive vertices is equal to that of negative vertices, hence we have
an even number of cycles.

On the twice punctured disk $D \cap X$, these edges $e_1, ..., e_n$
cuts it into one annulus and $n-1$ rectangles $D_1, ..., D_{2n-1}$.
As before, call the two components of $(V,K,r)|T$ the white region $W$
and the green region $G$, with the white region containing $\bdd V$.
Then the $n/2$ of rectangles $D_{2i-1}$ are in the green region.
Also, the Dehn filling solid torus $N(K_r)$ is cut by $T$ into $n$
components $H_1, ..., H_n$, with $H_{2i-1}$ in the green region, and
each $H_{2i-1}$ is incident to two of the rectangles.  It follows that
if we shrink each $H_{2i-1}$ to an arc $\alpha_i$ then $\cup (H_{2i-1}
\cup D_{2i-1})$ becomes a set of annuli or M\"obius bands containing
these $\alpha_i$, with boundary on the above cycle.  But since the two
ends of $\alpha_i$ are of opposite signs, which by the above are on
different cycles, we see that there is no M\"obius band in the
above union, so they are all annuli.

Let $A$ be one of these annuli.  Then $\bdd A$ cuts $T$ into two
annuli $A_1, A_2$.  Since we assumed that $(V,K,r)$ contains no Klein
bottle, each $A \cup A_i$ is a torus instead of Klein bottle, hence
the frontier of $N(T \cup A)$ consists of three tori $T_0 \cup T_1
\cup T_2$, with $T_1, T_2$ in the green region.  One can check that
each $T_i$ has fewer intersection with $K_r$ than $T$, hence by the
minimality of $n$ we see that each $T_i$ bounds a solid torus $V_i$,
which must be disjoint from $T \cup A$ as otherwise it would contain
$\bdd V$ and hence have at least two boundary components,
contradicting the assumption that $V_i$ are solid tori.  It now
follows that $G$ is homeomorphic to the manifold obtained by gluing
$V_1, V_2$ along an annulus.  The incompressibility of $T$ then
implies that $G$ is a small Seifert fiber space with orbifold a disk
with two cone points.  
\qed

\begin{thm} \label{thm:4.3}
  Suppose $\wrap(K) = 2$, and $(V,K,r)$ is non hyperbolic.  Then
  $K_n(r_n)$ is nonhyperbolic for all but at most three $n$.
  Moreover, either

  (1) there is an $n_0$ such that $K_n(r_n)$ is toroidal unless $|n -
  n_0| \leq 1$; or

  (2) $K_n(r_n)$ is atoroidal for all $n$, and there exist $q_1, q_2
  \in \Bbb Z$ such each $K_n(r_n)$ is either reducible or has a small
  Seifert fibration with $q_1, q_2$ as the indices of two of its
  singular fibers.
\end{thm}

\proof By \cite{Be,Sch} $(V,K,r)$ is irreducible and not a solid
torus, hence it is either a small Seifert fibered manifold or
toroidal.  If it is a small Seifert fibered manifold then (2) holds,
where $q_1, q_2$ are the indices of the two singular fibers of
$(V,K,r)$.

Suppose $(V,K,r)$ is toroidal.  Let $T$ be an essential torus of
$(V,K,r)$ given by Lemma \ref{lemma:4.2}.  Note that $K_n$ is obtained
by Dehn filling on one component of a hyperbolic link, hence by
\cite{Wu92} it is nontrivial for all but at most two adjacent integers
$n$.  By \cite{Ga3}, $K_n(r_n)$ cannot contain a nonseparating sphere
if $K_n$ is nontrivial.  Therefore if $T$ is nonseparating then it
remains a nonseparating incompressible torus in $K_n(r_n)$ for all but
at most two consecutive $n$, hence (1) holds.

We may now assume that $T$ is separating.  By Lemma \ref{lemma:4.2},
$T$ cuts $(V,K,r)$ into $M_1, M_2$, where $M_2$ contains $T_0 = \bdd
V$, and $M_1$ is a small Seifert fiber space.  Thus $K_n(r_n) = M_1
\cup_T M_2(r_n)$.  Let $q_1, q_2$ be the indices of the singular
fibers of $M_1$.  By \cite[Theorem 2.4.4]{CGLS}, if $M_2$ is not a cable
space then $T$ is incompressible in $M_2(r_n)$ and hence
incompressible in $K_n(r_n)$, for all but at most two consecutive
$n$, so again (1) follows and we are done.

We now assume that $M_2$ is a cable space.  Let $A$ be an essential
annulus in $M_2$ with one boundary component on each of $T$ and $T_0$,
and let $\gamma_0$ be the boundary slope of $A$ on $T_0$.  Let $(m,l)$
be a meridian-longitude pair of $T_0 = \bdd V$.  Then $K_n(r_n)$
is obtained from $(V,K,r)$ by Dehn filling on $T_0$ along the slope
$\alpha_n = l -nm$.

By \cite[Theorem 2.4.3]{CGLS}, there is a slope $\gamma_0$ on $T_0$
such that $T$ remains incompressible in $M_2(\alpha_n)$ if
$\Delta(\alpha_n, \gamma_0) \geq 2$.  If $m \neq \gamma_0$ then at
most three consecutive $\alpha_n$ satisfies the above condition and
hence (1) holds.  Now assume $m = \gamma_0$.  Then $M_2(\alpha_n)$ is
a solid torus for all $n$, so $K_n(r_n)$ is the union of the small
Seifert fiber space $M_1$ and the solid torus $M_2(\alpha_n)$.  If the
fiber of $M_1$ is the meridional slope of $M_2(\alpha_n)$ then
$K_n(r_n)$ is reducible, and if not then the Seifert fibration of
$M_1$ extends to a small Seifert fiber structure of $K_n(r_n)$.  Hence
(2) holds in this case.  \qed

\section{Surgery on wrapped Montesinos knots} \label{sec:5}

Denote by $T[t_1, ..., t_p]$ the Montesinos tangle consisting of $p$
rational tangles of slopes $t_i$; see Figure 5.1(a) for $p=2$, where a
circle with label $t_i$ represents a rational tangle of slope $t_i
\neq 1/0$.  Up to isotopy we may assume $t_i$ are not integers unless
$p=1$.  We can add two strings to connect the top endpoints to the
bottom ones to make it a knot of wrapping number 2 in a solid torus
$V$.  Up to homeomorphism of $V$ there are two ways to add these two
strings, as shown in Figure 5.1(b)-(c), denoted by $K^0(t_1, ...,
t_p)$ and $K^1(t_1, ..., t_p)$, respectively.  Recall that the circle
with a dot in these figures represents the component $K'$ to be
removed, so $V = S^3 - \Int N(K')$.  Only one of these is a knot if
the two top endpoints of the tangle belong to different strings.  We
call these knots {\em wrapped Montesinos knots\/} in solid tori.  Note
that if $K = K^a(t_1, ..., t_p)$ for $a =0, 1$ then $K_n$ is a
Montesinos knot $M(1/(a+2n), t_1, ..., t_p) = M(t_1, ..., t_p,
1/(a+2n))$ in $S^3$.  The purpose of this section is to determine all
exceptional Dehn surgeries on these wrapped Montesinos knots in solid
tori.

\bigskip
\leavevmode

\centerline{\epsfbox{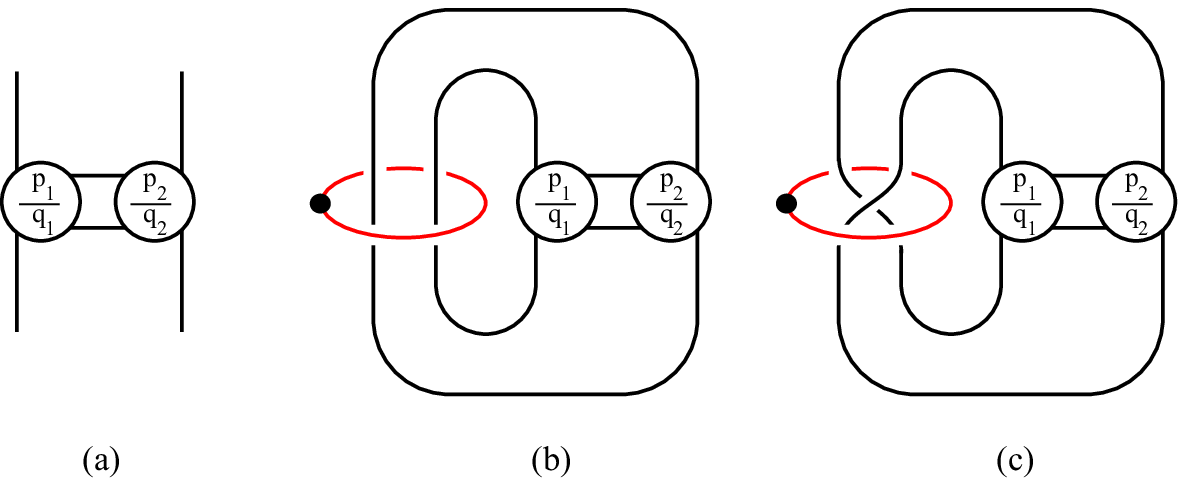}}
\bigskip
\centerline{Figure 5.1}
\bigskip

As before, we fix a meridian-longitude pair $(m_0,l_0)$ of $K \subset
V$ so that it becomes the preferred meridian-longitude pair of $K_0
\subset S^3$.  A slope $r$ is then represented by a rational number
$u/v$ if it represents $\pm (um_0 + vl_0)$ in $H_1(\bdd N(K))$.  We first
consider the knot $K = K^1(-1/2,\, 1/3)$.  By Proposition
\ref{prop:2.2}, $K(6)$ and $K(8)$ are toroidal, and $K(7)$ is small
Seifert fibered.  The following lemma shows that there is no other
exceptional surgery on this knot.

\begin{lemma} \label{lemma:5.1} Let $K = K^1(-1/2,\, 1/3)$.  Then
  $(V,K,r)$ is an exceptional surgery if and only if $r = 6,7,8$.
\end{lemma}

\proof Let $\varphi_n: V \to S^3$ be the embedding defined in Section
4.  Then $K_n = \varphi_n(K)$ is the $(-2,\, 3,\, 1+2n)$ pretzel knot,
and we have $r_n = \varphi_n(r) = r + 4n$ with respect to the
preferred meridian-longitude of $K_n$.  Assume that $(V,K,r)$ is
nonhyperbolic.  By Theorem \ref{thm:4.3} either (1) $K_n(r_n)$ is
toroidal for all but at most three $r_n$, or (2) $K_n(r_n)$ is
reducible or small Seifert fibered for all $n$.  If (1) is true then
by \cite{Wu:pp1} we have $r = 8$.  Hence we assume (2) holds.  When
$n=-1$, $K_n$ is the $(-2,3,-1)$ pretzel knot, which can be deformed
to the mirror image of the knot $5_2$ on the knot table.  By \cite{BW}
it admits only three Seifert fibered surgery, with slopes $1,2,3$,
respectively, and there is no reducible surgery. Since $r_n = r + 4n =
r - 4$ when $n=-1$, we have $r = 5, 6, 7$.

It remains to show that $(V,K,5)$ is hyperbolic.  Let $r=5$.  The
above shows that $K_n(r_n)$ is small Seifert fibered for $n=-1$.  For
$n=0, 1, 2$, the knot $K_n$ is the $(2,5)$, $(3,4)$ and $(3,5)$ torus
knot, respectively, hence $K_n(r_n)$ cannot be toroidal for these four
$n$.  By Theorem \ref{thm:4.3}, this implies that if $(V,K,5)$ is
nonhyperbolic then conclusion (2) of that theorem must hold, i.e.\
there exists $q_1, q_2$ such that each $K_n(r_n)$ is either reducible,
or has a small Seifert fibration with $q_1, q_2$ as the indices of
two of its singular fibers.

For $n=0$, $K_n$ is the $(2,5)$ torus knot, and $r_n = r = 5$.  The
cabling slope of $K_0$ is $2 \times 5=10$, hence the Seifert fiber
structure of the exterior of $K_0$ extends over the Dehn filling solid
torus, whose core is then a singular fiber of index $10-5 = 5$.
Therefore $K_0(5)$ is a small Seifert fibered manifold with three
singular fibers of indices $2, 5, 5$, respectively.  Similarly, $K_1$
is the $(-2,3,3)$ pretzel knot, which is the $(3,4)$ torus knot.  The
cabling slope of $K_1$ is $3 \times 4 = 12$, and the surgery slope is
$r_1 = 5+4 = 9$, so after Dehn surgery the manifold $K_1(r_1)$ is a
small Seifert fibered manifold with three singular fibers of indices
$3,4,3$, respectively.  By \cite[Theorem VI.17]{Ja}, Seifert
fibrations for these manifolds are unique.  This is then a
contradiction to Theorem \ref{thm:4.3} since no pair of indices of the
singular fibers of $K_0(r_0)$ match those of $K_1(r_1)$.  \qed

\begin{lemma} \label{lemma:5.2} Let $K = K^1(-1/2,\, 1/q)$, where $|q|
  \geq 3$ is odd.  Let $X$ be the tangle space as defined in Section
  \ref{sec:2}.  Then $F = \bdd_h X$ is incompressible unless $q = 3$.
\end{lemma}

\proof Let $\hat X$ be the manifold obtained by attaching a 2-handle
to $X$ along the annulus $\bdd_v X$, and let $\hat F = \bdd \hat X$ be
the corresponding surface.  Note that $X$ is a handlebody of genus 2,
so if $F$ is compressible then there is a nonseparating compressing
disk $D_1$, which remains a compressing disk in $\hat X$.  Let $L$ be
the link obtained by adding two horizontal arcs to the tangle $T[-1/2,
1/q]$.  Then $\hat X = E(L)$, hence $L$ is a trivial knot.  On the
other hand, it is easy to see that $L$ is a $(2,\, q-2)$ torus knot.
Since $|q| \geq 3$, it follows that $L$ is trivial if and only if
$q=3$.  
\qed 
\bigskip

The knot $K = K^a(1/q_1,\, 1/q_2)$ in $V$ has an obvious spanning
surface which is a once punctured torus or Klein bottle, called the
{\it pretzel surface}.  Its boundary slope is called the {\it pretzel
  slope} of $K$.

\begin{lemma} \label{lemma:5.3} Let $K = K^a(1/q_1,\, 1/q_2)$ be a
  pretzel knot in $V$, $|q_i|>1$ and $\{q_1,\, q_2\} \neq \{ \mp 2,\,
  \pm 3\}$.  Then $(V,K,r)$ is hyperbolic unless $r$ is the pretzel
  slope.
\end{lemma}

\proof By Lemma \ref{lemma:5.2} and \cite[Lemma 3.3]{Wu96} the surface
$\bdd_h X$ is incompressible, hence by Theorem \ref{thm:3.8} we see
that the knot $K \subset V$ admits no reducible or Seifert fibered
surgery and at most one toroidal surgery.  Since the surgery along the
pretzel slope $r$ contains either a nonseparating torus or a Klein
bottle and hence is nonhyperbolic, it is the only exceptional surgery
slope.  \qed

\begin{lemma} \label{lemma:5.4}
  Suppose $K = K^1(-1/2,\, 2/5)$.  Then $(V,K,r)$ admits no
  exceptional surgery.
\end{lemma}

\proof We have $K_n = M(-1/2,\, 2/5,\, 1/(1+2n))$.  Checking the list
in \cite[Theorem 1.1]{Wu:pp1}, we see that $K_n$ admits no toroidal
surgery when $n>9$, so by Theorem \ref{thm:4.3}, if $(V,K,r)$ is an
exceptional surgery then $K_n(r_n)$ is atoroidal and nonhyperbolic for
all $n$.  In particular, this should be true for $n=-1$, in which case
$K_n = M(-1/2,\, 2/5,\, -1)$ can be deformed to the 2-bridge knot
associate to the rational number $-1/(3-1/4) = -4/11$.  On the other
hand, by \cite{BW} this knot admits only one exceptional surgery,
which produces a toroidal manifold.  Hence we have a contradiction.
\qed \bigskip

Two Montesinos tangles $T[t_1, ..., t_k]$ and $T[s_1, ..., s_k]$ are
equivalent if $s_i - t_i$ are integers, and $\sum s_i = \sum t_i$, in
which case $K^a(t_1, ..., t_k)$ is isotopic to $K^a(s_1, ..., s_k)$.
Any $t_i = 0$ can be added or deleted without affecting the knot type.
Note that $K = K^a(t_1, ..., t_k)$ is isotopic to $K' = K^a(t_k, ...,
t_1)$, and is the mirror image of $K'' = K^a(-t_1, ..., -t_k)$, so
$(V,K,r)$ is homeomorphic to $(V,K',r)$ and $(V,K'',-r)$.  When $k=1$,
twisting $m$ times along a meridional disk of $V$ will change $K =
K^a(t_1)$ to $K''' = K^a(1/(2m + 1/t_1))$.  We will consider these
knots $K, K', K'', K'''$ as {\it equivalent}.  We may assume that $K$
is not equivalent to $K^a(0)$ or $K^a(1/q)$ as otherwise $K$ is
nonhyperbolic.  The following theorem classifies exceptional surgeries
on wrapped Montesinos knots.

\begin{thm} \label{thm:5.5}
  Suppose $K = K^a(t_1, ..., t_k) \subset V$ is not equivalent to
  $K^a(0)$ or $K^a(1/q)$ for any integer $q$.  Let $(V,K,r)$ be an
  exceptional surgery.  Then $(K, r)$ is equivalent to one of the
  following pairs.  The surgery is small Seifert fibered for $r=1,2,3$
  in (1) and $r=7$ in (4), and toroidal otherwise.

  (1) $K = K^0(2)$ (the Whitehead knot), $r = 0, 1, 2, 3, 4$.

  (2) $K = K^a(n)$, $n>2$, $r=0$ if $a=0$, and $r=2n$ otherwise.

  (3) $K = K^a(1/q_1,\, 1/q_2)$, $|q_i|\geq 2$, and $r$ is the pretzel
  slope.

  (4) $K = K^1(-1/2,\, 1/3)$, $r = 6,7,8$.
\end{thm}

\proof First assume that $k = 1$, so $T[t]$ is a rational tangle.  In
this case any $K^1(t')$ is equivalent to some $K^0(t)$.  By the above,
the reciprocal $1/t$ has the property that $K^0(t)$ is equivalent to
$K^0(t')$ if $1/t' = 2 \pm (1/t)$, and by assumption $1/t \neq 0, 1$.
Hence up to equivalence we may assume that $0 < 1/t < 1$, i.e.\ $t =
p/q > 1$.  Note that $K_n$ is the 2-bridge knot in $S^3$ associated to
the rational number $r = 1/(2n + q/p)$.  Hence if $q \neq 1$ then for
any $n>1$, $K_n$ is not equivalent to a 2-bridge knot associated to
any rational number of type $1/(b_1 + 1/b_2)$ with $b_1, b_2 \in \Bbb
Z$.  It follows from \cite[Theorem 1.1]{BW} that all nontrivial
surgeries on such $K_n$ are hyperbolic.  By Theorem \ref{thm:4.3} this
implies that $K \subset V$ admits no exceptional surgery.  For $q=1$,
$t = p/q > 1$ is an integer.  If $t = 2$ then $K$ is the Whitehead
knot in $V$ and it is well known that $K$ admits exactly 5 exceptional
surgeries as listed in (1).  The hyperbolicity of $(V,K,r)$ for $r\neq
0, ..., 4$ can also be proved using the argument in the proof of Lemma
5.4 and the classification of exceptional surgeries on 2-bridge knots
given in \cite{BW}.  If $t>2$ then $K$ is not the Whitehead knot, and
the argument in the proof of Lemma \ref{lemma:5.3} shows that $\bdd_h
X$ is incompressible, hence by Theorem \ref{thm:3.8} we see that $K$
admits no Seifert fibered surgery and at most one toroidal surgery.
Note that $K$ has a spanning surface $F$ in $V$ which is a once
punctured torus or Klein bottle.  As in the proof of Lemma 5.3,
surgery long the boundary slope of $F$ produces a toroidal manifold,
so there is no other exceptional surgery.  The toroidal slope is given
in (2).

We now consider the case that $k > 1$.  We may assume that $q_i \geq
2$ for all $i$ as otherwise the Montesinos tangle is equivalent to one
with fewer rational tangles.  If $k\geq 3$ then $K_n$ is a Montesinos
know of length at least 4 for all $|n|\geq 2$.  By \cite{Wu96jdg}
$K_n$ admits no exceptional surgery.  Hence by Theorem \ref{thm:4.3}
we see that $K(r)$ is hyperbolic for all nontrivial $r$, so there is
no exceptional surgery.

We now assume $K = K^a(p_1/q_1,\, p_2/q_2)$ with $q_i \geq 2$.  Then
$K_n = K(p_1/q_1,\, p_2/q_2,\, 1/2n)$ or $K(p_1/q_1,\, p_2/q_2,\,
1/(2n+1))$.  By Theorem \ref{thm:4.3}, if $(V,K,r)$ is exceptional
then either $K_n(r_n)$ is toroidal for all but at most three $n$, or
it is either reducible or atoroidal and Seifert fibered for all $n$.  By
\cite{Wu:pp1}, if $K(p_1/q_1,\, p_2/q_2,\, 1/q_3)$ admits a toroidal surgery
and $|q_3| > 9$ then $|p_i| = 1$ and the surgery slope is the pretzel
slope.  Hence if $K_n(r_n)$ is toroidal for almost all $n$ then $K =
K^a(1/q_1,\, 1/q_2)$ and $r$ is the pretzel slope, so (3) holds.

We may now assume that $K_n(r_n)$ is reducible or atoroidal and
Seifert fibered for all $n$.  As above, we have $K_n =
M(p_1/q_1,\, p_2/q_2,\, 1/2n)$ or $K_n = M(p_1/q_1,\, p_2/q_2,\,
1/(2n+1))$, and by \cite{Wu96jdg} $K_n(r_n)$ cannot be reducible; hence it
must be an atoroidal small Seifert fibered manifold for any $n$.  By
\cite[Theorems 7.2 and 7.3]{Wu:pp2}, one of the following must hold.

(i) $K_n$ is a $(q_1, q_2, q_3, d)$ pretzel knot or its mirror image,
and either $d=0$, or all $q_i$ are positive and $d=-1$.  Moreover,
either some $|q_i| = 2$ or $|q_i| = |q_j| = 3$ for some $i\neq j$.

(ii) $K_n = K(\mp 1/2,\, \pm 2/5,\, 1/(2n+1))$. 

In (i) above, the case $d=-1$ cannot happen in our case since
$K_n(r_n)$ is atoroidal Seifert fibered for both $n$ positive and
negative, contradicting the condition that all $q_i$ are of the same
sign (up to taking mirror image of $K_n$.)  Therefore $K_n$ must be a
genuine pretzel knot if (i) holds.  It follows that the tangle
must be equivalent to $T[1/q_1, 1/q_2]$ or $T[-1/2, 2/5]$.

The result now follows from Lemmas \ref{lemma:5.3} and \ref{lemma:5.4}.
\qed

\bigskip

\noindent
Department of Mathematics,  University of Iowa,  Iowa City, IA 52242
\\
Email: {\it wu@math.uiowa.edu}

\enddocument